\documentclass[12pt]{amsart}
\usepackage{amsfonts}
\newtheorem{theorem}{Theorem}[section]
\newtheorem{lemma}[theorem]{Lemma}
\newtheorem{proposition}[theorem]{Proposition}
\newtheorem{corollary}[theorem]{Corollary}
\newtheorem{definition}[theorem]{Definition}
\newtheorem{example}[theorem]{Example}

\newtheorem{remark}[theorem]{Remark}

\begin{document}

\title{On Pairwise paralindeloff Bitopological Spaces}

\author{Hend Bouseliana}
 \address{Department of Mathematics, University Putra Malaysia, 43400 UPM, Serdang, Selangor
 Malaysia, E-mail:bouseliana2007@yahoo.com}

\author{Adem K\i l\i \c{c}man}
 \address{Department of Mathematics and Institute for Mathematics Research,
  University Putra Malaysia, 43400 UPM, Serdang, Selangor, Malaysia, E-mail:akilic@upm.edu.my}

\begin{abstract} The aim of this work is to introduce and study some new types of generalizations of pairwise paralindeloff spaces, pairwise nearly paralindeloff and almost paralindeloff spaces. Some of their characterizations, properties and subsets are investigated and studied. Furthermore, the relation of among them with some pairwise separation axioms are considered.\\

\textbf{Keywords and phrases}: Bitopological spaces, pairwise regular open sets, separation axioms, paralindeloff spaces.\\

\textbf{2010 AMS Subject Classification:} 54A05; 54D10; 54C20 \end{abstract}
\maketitle

\section{Introduction and Preliminaries}

\noindent The idea of bitopological spaces was first initiated by Kelly \cite{Kelly63} then thereafter several studies were carried out in order to generalize the topological concepts to bitopological settings. The concept of nearly paralindeloff spaces was first introduced and studied by Thanapalan \cite{DT}. Further, in \cite{DT92} and \cite{DT94}, studied the idea of almost paralindeloff space and presented some of properties. The aim of this work is to study the idea of nearly and almost paralindeloff bitopological spaces, which will be denoted by $B_1$-nearly paralindeloff, $B_r$-nearly paralindeloff, $B_1$-almost paralindeloff spaces. Some of their characterizations, properties, subspaces and subsets are studied and investigated. In addition, the relation of among them and some pairwise separation axioms are also considered.\\

\noindent Throughout this paper, all spaces $(X,\tau)$ and $(X,\tau_1,\tau_2)$ are always meant topological spaces and bitopological spaces, respectively. A space $(X,\tau_1,\tau)_2)$ is said to be $(i,j)$-almost regular \cite{SS} (resp. $(i,j)$-regular space \cite{Kelly63}) if for each $x \in X$ and each $(i,j)$-regular open (resp. $i$-open) set $U$ containing $x$, there exists an $(i,j)$-regular open (resp. $i$-open) set $V$ such that $x \in V \subseteq j-cl(V) \subseteq U$. $(X,\tau_1,\tau)_2)$ is called pairwise almost regular (resp. pairwise regular) if it is $(1,2)$-almost regular and $(2,1)$-almost regular (resp. $(1,2)$-regular and $(2,1)$-regular).

\begin{definition} A space $(X,\tau_1,\tau_2)$ is said to be

\begin{enumerate}

\item $i$-$P$-space if any countable intersection of arbitrary collection of $i$-open sets in $X$ is $i$-open. $X$ is called $P$-space if it is $i$-$P$-space for each $i=1,2$.

\item $(i,j)$-$P$-space if any countable intersection of arbitrary collection of $i$-open sets in $X$ is $j$-open. $X$ is called $B$-$P$-space if it is both $(1,2)$-$P$-space and $(2,1)$-$P$-space.

\item $(i,j)$-$P_r$-space if any countable intersection of arbitrary collection of $(i,j)$-regular open sets in $X$ is $(i,j)$-regular open. $X$ is called pairwise $P_r$-space if it is both $(1,2)$-$P_r$-space and $(2,1)$-$P_r$-space.

\item $(i,j)$-weakly $P$-space if for any countable family $\{U_n:n \in \mathbb{N}\}$ of $i$-open sets in $X$, we have
$j-cl(\bigcup_{n \in \mathbb{N}}) U_n= \bigcup_{n \in \mathbb{N}} j-cl(U_n)$. $X$ is called pairwise weakly $P$-space if it is both $(1,2)$-weakly $P$-space and $(2,1)$-weakly $P$-space.
\end{enumerate} \end{definition}

\begin{definition} \cite{Zabidin1}\cite{Zabidin6}. A bitopological space $X$ is said to be $(i,j)$-nearly Lindelof if for every $i$-open cover $\{U_\alpha :\alpha \in \Delta\}$ of $X$, there exists a countable subset $\{\alpha_n:n \in \mathbb{N}\}$ such that $X= \bigcup_{n \in \mathbb{N}} i-int(j-cl(U_{\alpha_n}))$. $X$ is called pairwise nearly Lindelof if it is both $(1,2)$-nearly Lindelof and $(2,1)$-nearly Lindelof.\end{definition}

The concept of extremely disconnectedness was first extended to a bitopological space $(X,\tau_1,\tau_2)$ by Datta \cite{Datta1}.

\begin{definition} \cite{Datta1}. A bitopological space $(X,\tau_1,\tau_2)$ is said to be $(i,j)$-extremely disconnected if the $i$-closure of every $j$-open set is $j$-open. $X$ is called pairwise extremely disconnected if it is both $(1,2)$-extremely disconnected and $(2,1)$-extremely disconnected.\end{definition}

\section{pairwise nearly paralindelof spaces}

\begin{definition} A bitopological space $X$ is said to be $(i,j)_1$-nearly paralindelof if every cover of $X$ by $(i,j)$-regular open sets has $i$-open refinement which is $j$-locally countable. $X$ is called $B_1$-nearly paralindelof if it is both $(1,2)_1$-nearly paralindelof and $(2,1)_1$-nearly paralindelof. \end{definition}

\begin{remark} It is clear that every $(i,j)$-nearly Lindelof is $(i,j)_1$-nearly paralindelof.\end{remark}

\begin{example} Let $\mathbb{Q^c}$ be an irrational numbers with the discrete topology $\tau_{dis}$. It is clear that  $(\mathbb{Q^c},\tau_{dis},\tau_{dis})$ is $B_1$-nearly paralindelof (see \cite{SA}). But $X$ is not pairwise nearly Lindelof space because $X$ is not pairwise almost Lindelof (see \cite{Zabidin6}).\end{example}

Before the next theorem, we state the some properties of $i$-locally countable of a bitopological space $X$.

\begin{lemma}  \label{l1} Let $(X,\tau_1,\tau_2)$ be a bitopological space and $\mathcal{B}$ be a basis containing $\tau_i$-open subsets of $X$. Then $A$ is $i$-open if and only if for each point $x \in A$ there is $ U \in \mathcal{B}$ with $x \in U \subset A$. \end{lemma}

\noindent In particular, since $\mathcal{B}$ a basis for $\tau_i$, Lemma \ref{l1} provides a very useful way for showing that a given set $A$ is $i$-open for $i=1,2.$\\

\noindent The main property of $i$-locally countable families is following.\\

\begin{theorem}  \label{locally countable} Let a bitopological space $(X,\tau_1,\tau_2)$ be $i-P$-space and $\mathcal{U}=\{U_\alpha| \alpha \in \Delta\}$ be $i$-locally countable family in $X$. Then:
\begin{enumerate}
\item $\{i-cl(U_\alpha): \alpha \in \Delta\}$ is also $i$-locally countable.
\item For each $\Delta^{'} \subseteq \Delta$, $\bigcup \{i-cl(U_\alpha): \alpha \in \Delta^{'}\}$
 is $i$-closed in $X$.
\end{enumerate}\end{theorem}

\textbf{Proof} (1). Since the collection $\{U_\alpha: \alpha \in \Delta\}$ is $i$-locally countable, then for each $x \in X$ there is a neighborhood $O_x$ such that $U_\alpha \bigcap O_x=\emptyset$ for all but at most countably many $\alpha$.
So $U_\alpha \bigcap O_x=\emptyset \Rightarrow U_\alpha \subseteq X-O_x \Rightarrow i-cl(U_\alpha) \subseteq X-O_x \Rightarrow i-cl(U_\alpha) \cap O_x=\emptyset$. This completes the proof.\\

(2). Let $B = \bigcup_{\alpha \in \Delta} i-cl(U_\alpha)$. From (1), for each $x \notin B$ there is a neighborhood $O_x$ meeting at most countably many members of $\{i-cl(U_\alpha): \alpha \in \Delta\}$ say, $i-cl(U_{\alpha_1}), i-cl(U_{\alpha_2}),..., i-cl(U_{\alpha_n}),...$; then $O_x \bigcap (\bigcap_{n \in \mathbb{N}} (X-i$-$cl(U_{\alpha_n})))$ is $i$-open neighborhood of $x$ not meeting $B$, so by Lemma \ref{l1}, $X-B$ is $i$-open. Thus $B$ is $i$-closed set.

\begin{corollary} Let a bitopological space $(X,\tau_1,\tau_2)$ be $P$-space and $\mathcal{U}=\{U_\alpha: \alpha \in \Delta\}$ be locally countable family in $X$. Then:
\begin{enumerate}
\item $\{cl(U_\alpha)| \alpha \in \Delta\}$ is also locally countable.
\item For each $\Delta^{'} \subset \Delta$, $\bigcup \{cl(U_\alpha): \alpha \in \Delta^{'}\}$
 is closed in $X$.
\end{enumerate} \end{corollary}

\begin{theorem} \label{th4.1} If $X$ is $(i,j)_1$-nearly paralindelof, $(i,j)$-almost regular and $j$-$P$-space, then:

\begin{enumerate}
\item Every $(i,j)$-regular open cover of $X$ has $(i,j)$-regular open refinement which is $j$-locally countable.
\item Every $(i,j)$-regular open cover of $X$ has $j$-closed refinement which is $j$-locally countable.
\item Every $(i,j)$-regular open cover of $X$ has $(i,j)$-regular closed refinement which is $j$-locally countable.
\end{enumerate}
\end{theorem}
\textbf{Proof.}

$(1)$. Let $\mathcal{U}$ be any $(i,j)$-regular open cover of $X$. Let $x \in X$. Then, there exists some $U \in \mathcal{U}$ containing $x$. Since $X$ is $(i,j)$-almost regular, there is $(i,j)$-regular open set $V_x$ such that $x \in V_x \subseteq j-cl(V_x) \subseteq U_\alpha$. Now, $\mathcal{V}=\{V_x :x \in X \}$ is $(i,j)$-regular open cover of $X$ therefore it has $j$-locally countable family $\mathcal{W}=\{W_\lambda: \lambda \in \Lambda\}$ of $i$-open sets in $X$ refining $\mathcal{V}$ and covering $X$. It is known that $W_\alpha \subseteq i-int(j-cl(W_\alpha)) \subseteq j-cl(W_\alpha)$. Set $\mathcal{H}=\{i-int(j-cl(W_\alpha)): \alpha \in \Delta\}$. So, $\mathcal{H}$ covers $X$ and $j$-locally countable because it is a subfamily of $j$-locally countable by Theorem \ref{locally countable}. Moreover,  since $\mathcal{W}$ is refinement of $\mathcal{V}$, then for each $W_\alpha \in \mathcal{W}$ there is a $V_{x(\alpha)} \in \mathcal{V}$ such that $W_\alpha \subseteq V_{x(\alpha)}$. Furthermore, we have $i-int(j-cl(W_\alpha)) \subseteq j-cl(W_\alpha) \subseteq j-cl(V_{x(\alpha)}) \subset U$ for some $U \in \mathcal{U}$. Therefore, $\mathcal{H}$ is also refinement of $\mathcal{U}$. This completes the proof.\\

 $(2)$. Let $\mathcal{U}$ be any $(i,j)$-regular open cover of $X$. Let $x \in X$. Then, there exists some $U \in \mathcal{U}$ containing $x$. Since $X$ is $(i,j)$-almost regular, there is $(i,j)$-regular open set $V_x$ such that $x \in V_x \subseteq j-cl(V_x) \subseteq U_\alpha$. Now, $\mathcal{V}=\{V_x :x \in X \}$ is $(i,j)$-regular open cover of $X$ therefore it has $j$-locally countable family $\mathcal{W}=\{W_\lambda: \lambda \in \Lambda\}$ of $i$-open sets in $X$ refining $\mathcal{V}$ and covering $X$. Let $\mathcal{H}=\{j-cl(W):\alpha \in \Delta\}$ so that is also $j$-locally countable by Theorem \ref{locally countable}. Since $\mathcal{W}$ is refinement of $\mathcal{V}$, then for each $W_\alpha \in \mathcal{W}$ there is a $V_{x(\alpha)} \in \mathcal{V}$ such that $W_\alpha \subseteq j-cl(W_\alpha) \subset j-cl(V_{x(\alpha)}) \subseteq U$ for some $U \in \mathcal{U}$. Therefore, $\mathcal{H}$ is refinement of $\mathcal{U}$ and covers $X$ and the proof finishes.\\

$(3)$. Let $\mathcal{U}$ be any $(i,j)$-regular open cover of $X$. Let $x \in X$. Then, there exists some $U \in \mathcal{U}$ containing $x$. Since $X$ is $(i,j)$-almost regular, there is $(i,j)$-regular open set $V_x$ such that $x \in V_x \subseteq j-cl(V_x) \subseteq U_\alpha$. Now, $\mathcal{V}=\{V_x :x \in X \}$ is $(i,j)$-regular open cover of $X$ therefore it has $j$-locally countable family $\mathcal{W}=\{W_\lambda: \lambda \in \Lambda\}$ of $i$-open sets in $X$ refining $\mathcal{V}$ and covering $X$. Now, since $\mathcal{W}$ is a collection of $i$-open sets, then we can rewrite $\mathcal{H}$ in $(2)$ as following $\mathcal{H}=\{j-cl(W_\alpha): \alpha \in \Delta\}=\{j-cl(i-int(W_\alpha)): \alpha \in \Delta\}$ which is a family of $(j,i)$-regular closed sets, $j$-locally countable refining $\mathcal{U}$ and covering $X$ as desire.\\

\begin{corollary} If $X$ is $B_1$-nearly paralindelof, pairwise almost regular and $P$-space, then:

\begin{enumerate}
\item Every pairwise regular open cover of $X$ has pairwise regular open refinement which is locally countable.
\item Every pairwise regular open cover of $X$ has closed refinement which is locally countable.
\item Every pairwise regular open cover of $X$ has pairwise regular closed refinement which is locally countable.
\end{enumerate}
\end{corollary}

The Theorem \ref{th4.1} leads to theorem and two corollaries as following.

\begin{theorem} In $(i,j)$-almost regular and $j$-$P$-space $X$, $X$ is $(i,j)_1$-nearly paralindelof if and only if every $(i,j)$-regular open cover of $X$ has $(i,j)$-regular open refinement which is $j$-locally countable.\end{theorem}

\textbf{Proof.} Let $X$ be $(i,j)_1$-nearly paralindelof. Then, the condition holds from Theorem \ref{th4.1}. Conversely, since each $(i,j)$-regular open set is $i$-open, the proof completes.\\

\begin{corollary} In pairwise almost regular and $j$-$P$-space $X$, $X$ is $B_1$-nearly paralindelof if and only if every pairwise regular open cover of $X$ has pairwise regular open refinement which is $j$-locally countable for $i,j=1,2$ and $i \neq j$.\end{corollary}

\begin{corollary} Let a bitopological space $X$ be $(i,j)$-almost regular and $j$-$P$-space. If every $(i,j)$-regular open cover of $X$ has $(i,j)$-regular open refinement which is $j$-locally countable, then every $(i,j)$-regular open cover of $X$ has $j$-closed refinement which is $j$-locally countable.\end{corollary}

\textbf{Proof.} Since $i-int(j-cl(A)) \subseteq j-cl(A)$ and by using Theorem \ref{locally countable}, then with the same process of Theorem \ref{th4.1} $(1)$, the  proof is complete.\\

\begin{corollary} Let a bitopological space $X$ be pairwise almost regular and $P$-space. If every pairwise regular open cover of $X$ has pairwise regular open refinement which is locally countable, then every pairwise regular open cover of $X$ has closed refinement which is locally countable.\end{corollary}

\begin{corollary} Let a bitopological space $X$ be $(i,j)$-almost regular and $j$-$P$-space. If every $(i,j)$-regular open cover of $X$ has $j$-closed refinement which is $j$-locally countable and covers $X$, then every $(i,j)$-regular open cover of $X$ has $(j,i)$-regular closed refinement which is $j$-locally countable (not necessary covering $X$).\end{corollary}

\begin{corollary} Let a bitopological space $X$ be pairwise almost regular and $P$-space. If every pairwise regular open cover of $X$ has closed refinement which is locally countable, then every pairwise regular open cover of $X$ has pairwise regular closed refinement which is locally countable (not necessary covering $X$).\end{corollary}

\begin{lemma} \label{semiregular} Let $(X,\tau_1,\tau_2)$ be a bitopological space and let $(X,\tau^{s}_{(1,2)},\tau^{s}_{(2,1)})$ its pairwise semiregularization. Then

\begin{enumerate}
	\item $\tau_i-int(C)=\tau^{s}_{(i,j)}-int(C)$ for every $\tau_j$-closed set $C$;
	\item $\tau_i-cl(A)=\tau^{s}_{(i,j)}-cl(A)$ for every $A \in \tau_j$;
	\item the family of $(\tau_i,\tau_j)$-regular open sets of $(X,\tau_1,\tau_2)$ is the same as the family of $(\tau^{s}_{(i,j)},\tau^{s}_{(j,i)})$-regular open sets of $(X,\tau^{s}_{(1,2)},\tau^{s}_{(2,1)})$;
	\item the family of $(\tau_i,\tau_j)$-regular closed sets of $(X,\tau_1,\tau_2)$ is the same as the family of $(\tau^{s}_{(i,j)},\tau^{s}_{(j,i)})$-regular closed sets of $(X,\tau^{s}_{(1,2)},\tau^{s}_{(2,1)})$;
	\item $(\tau^{s}_{(i,j)})^{s}_{(i,j)}=\tau^{s}_{(i,j)}$.
\end{enumerate} \end{lemma}

\begin{lemma} \label{2} If the family $\{V_\lambda :\lambda \in \Lambda\}$ is $i$-open refinement $j$-locally countable of $(i,j)$-regular open subsets $\{U_\alpha:\alpha \in \Delta\}$ in $j$-$P$-space $(X,\tau_1,\tau_2)$, then $\{i-int(j-cl(V_\lambda)): \lambda \in \Lambda\}$ is $i$-open $j$-locally countable of $\{U_\alpha:\alpha \in \Delta\}$ in $(X,\tau^s_{(1,2)},\tau^s_{(2,1)})$.\end{lemma}

\textbf{Proof:} Since $i-int(j-cl(V_\lambda)) \in \tau^s_{(i,j)}$ for every $\lambda \in \Lambda$ and $V_\lambda \subseteq i-int(j-cl(V_\lambda)) \subseteq i-int(j-cl(U_\alpha))=U_\alpha$. Then, $\{i-int(j-cl(V_\lambda)):\lambda \in \Lambda\}$ is $i$-open refinement of $\{U_\alpha :\alpha \in \Delta\}$ in $(X,\tau^s_{(1,2)},\tau^s_{(2,1)})$. Now, we need to show that $\{i-int(j-cl(V_\lambda)):\lambda \in \Lambda\}$ is $j$-locally countable. Let $x \in X$ and $U_x \in \tau^s_{(j,i)}$ containing $x$. Since $\tau^s_{(j,i)} \subseteq \tau_j$, $U_x \in \tau_j$. Now, $\{\lambda \in \Lambda: U_x \cap V_\lambda \neq \emptyset\}$ is countable because $\mathcal{V}$ is $j$-locally countable. Then, $\{\lambda \in \Lambda: U_x \cap j-cl(V_\lambda) \neq \emptyset\}$ is also $j$-locally countable by Theorem \ref{locally countable}. But,

\begin{center}

$\{\lambda \in \Lambda: U_x \cap i-int(j-cl(V_\lambda)) \neq \emptyset\} \subseteq \{\lambda \in \Lambda: U_x \cap j-cl(V_\lambda) \neq \emptyset\}$,
\end{center}
then $\{\lambda \in \Lambda: U_x \cap i-int(j-cl(V_\lambda)) \neq \emptyset\}$ is countable. Hence, $\{i-int(j-cl(V_\lambda)): \lambda \in \Lambda\}$ is $j$-locally countable family in $(X,\tau^s_{(1,2)},\tau^s_{(2,1)})$. This implies that $\{i-int(j-cl(V_\lambda)): \lambda \in \Lambda\}$ is $j$-locally countable $i$-open refinement of $\mathcal{U}$ in $(X,\tau^s_{(1,2)},\tau^s_{(2,1)})$. This completes our proof of Lemma \ref{2}.

\begin{corollary} If the family $\{V_\lambda :\lambda \in \Lambda\}$ is open refinement locally countable of pairwise regular open subsets $\{U_\alpha:\alpha \in \Delta\}$ in $P$-space $(X,\tau_1,\tau_2)$, then $\{i-int(j-cl(V_\lambda)): \lambda \in \Lambda\}$ is locally countable open refinement of $\{U_\alpha:\alpha \in \Delta\}$ in $(X,\tau^s_{(1,2)},\tau^s_{(2,1)})$ for $i,j=1,2$ and $i \neq j$.\end{corollary}

\begin{theorem} Let $(X,\tau_1,\tau_2)$ be $j$-$P$-space. Then $(X,\tau_1,\tau_2)$ is $(i,j)_1$-nearly paralindelof if and only if $(X,\tau^s_{(1,2)},\tau^s_{(2,1)})$ is $(i,j)_1$-nearly paralindelof.\end{theorem}

\textbf{Proof.} Assume that $(X,\tau_1,\tau_2)$ is $(i,j)_1$-nearly paralindelof and let $\mathcal{U}=\{U_\alpha: \alpha \in \Delta\}$ be a $(i,j)$-regular open cover of $(X,\tau^s_{(1,2)},\tau^s_{(2,1)})$. Then, using Lemma \ref{semiregular}, we get $\mathcal{U}$ is also $(i,j)$-regular open cover of $(X,\tau_1,\tau_2)$, hence it admits $j$-locally countable family $\mathcal{V}=\{V_\lambda : \lambda \in \Lambda\}$ of $i$-open sets in $(X,\tau_1,\tau_2)$ which refines $\mathcal{U}$. Since $V_\lambda \subseteq i-int(j-cl(V_\lambda))$ and using Lemma \ref{2}, the collection $\{i-int(j-cl(V_\lambda)):\lambda \in \Lambda\}$ is still $i$-open $j$-locally countable refinement of $\mathcal{U}$ in $(X,\tau^s_{(1,2)},\tau^s_{(2,1)})$. This leads to that $(X,\tau^s_{(1,2)},\tau^s_{(2,1)})$ is $(i.j)_1$-nearly paralindelof.\\
Conversely, suppose that $(X,\tau^s_{(1,2)},\tau^s_{(2,1)})$ is $(i.j)_1$nearly paralindelof and let $\mathcal{U}=\{U_\alpha \alpha \in \Delta\}$ be a $(i,j)$-regular open cover of $(X,\tau_1,\tau_2)$. So, $U_\alpha \in \tau^s_{(i,j)}$ for all $\alpha \in \Delta$. Then, by Lemma \ref{semiregular}, $\mathcal{U}$ is also $(i,j)$-regular open cover of $(X,\tau^s_{(1,2)},\tau^s_{(2,1)})$. Then, it has $i$-open $j$-locally countable refinement $\mathcal{V}=\{V_\lambda : \lambda \in \Lambda\}$ in $(X,\tau^s_{(1,2)},\tau^s_{(2,1)})$. However, $\tau^s_{(i,j)} \subseteq \tau_i$ implies that $\mathcal{V}$ is $i$-open refinement of $\mathcal{U}$ in $(X,\tau_1,\tau_2)$. To prove that $\mathcal{V}$ is $j$-locally countable in $(X,\tau_1,\tau_2)$, let $x \in X$ and $U_x$ be any $j$-open neighborhood of $x$ in $(X,\tau_1,\tau_2)$. Since $U_x \subseteq j-int(i-cl(U_x))$, $\{\lambda \in \Lambda:U_x \cap V_\lambda \neq \emptyset\} \subseteq \{\lambda \in \Lambda:j-int(i-cl(U_x)) \cap V_\lambda \neq \emptyset\}$. Since $j-int(i-cl(U_x)) \in \tau^s_{(j,i)} \subseteq \tau_j$ and $\mathcal{V}$ is $j$-locally countable in $(X,\tau^s_{(1,2)},\tau^s_{(2,1)})$, $\{\lambda \in \Lambda:j-int(i-cl(U_x)) \cap V_\lambda \neq \emptyset\}$ is countable. Thus, $\{\lambda \in \Lambda:U_x \cap V_\lambda \neq \emptyset\}$ is countable. therefore, $\mathcal{V}$ is $j$-locally countable in $(X,\tau_1,\tau_2)$. Then, $(X,\tau_1,\tau_2)$ is $(i,j)_1$-nearly paralindelof.\\

\begin{corollary} Let $(X,\tau_1,\tau_2)$ be $P$-space. Then $(X,\tau_1,\tau_2)$ is $B_1$-nearly paralindelof if and only if $(X,\tau^s_{(1,2)},\tau^s_{(2,1)})$ is $B_1$nearly paralindelof.\end{corollary}

\begin{theorem} Let $X$ be $(i,j)$-almost regular, $(i,j)_1$-nearly paralindelof, $(i,j)$-$P_r$-space, $(j,i)$-$P$ and $j$-$P$-space. Then, every $(i,j) $-regular open covering of $X$ has $i$-open star refinement.\end{theorem}

\textbf{Proof.} Let $\mathcal{U}=\{U_\alpha: \alpha \in \Lambda\}$ be any $(i,j)$-regular open covering of $X$. Then by Theorem \ref{th4.1}, $\mathcal{U}$ admits $j$-locally countable family $\mathcal{V}$ of $j$-closed refinement of $\mathcal{U}$. For each $V \in \mathcal{V}$, there exists a $U_V \in \mathcal{U}$ such that $V \subseteq U_V$. For each $x \in X$, let

\begin{center}
$M_x=\bigcap \{U_v:x \in V \in \mathcal{V}\} \bigcap \{\bigcap \{X-V:x \notin V \in \mathcal{V}\}\}$.

\end{center}

 Since $X$ is $(j,i)$-$P$-space, $\{\bigcap \{X-V:x \notin V \in \mathcal{V}\}\}$ is $i$-open set. Moreover, $H=\bigcap \{U_v:x \in V \in \mathcal{V}\}$ is $(i,j)$-regular open set because $X$ is $(i,j)$-$P_r$-space so that $H$ is $i$-open set. Thus, each $M_x$ is $i$-open set. Let define the family $\mathcal{M}=\{M_x:x \in X \}$. For $x \in X$, either $x \in V \in \mathcal{V}$ or $x \notin V \mathcal{V}$. So, in both cases $x \in M_x$ from that the family $\mathcal{M}$ forms $i$-open covering of $X$. Now we shall show that $\mathcal{M}$ is star refinement of $\mathcal{U}$. Suppose that $p \in X$ is a fixed point and $\mathcal{M}_p=\bigcup \{M_x :p \in M_x \in \mathcal{M}\}$. Since $\mathcal{V}$ covers $X$, there exists $V \in \mathcal{V}$ such that $p \in V$. We claim that $p \in M_x \subseteq U_V$. Let $M_x$ be any member of $\mathcal{M}$ containing $p$. By definition of $M_x$, if $x \notin V$, then $p \in M_x \subseteq X-V$. But this contradicts $p \in V$. Therefore, $x \in V$ and $M_x \subseteq U_V$. Then, every member of $\mathcal{M}$ containing $p$ is contained in $U_V$ and $\mathcal{M}_p \subseteq U_V \in \mathcal{U}$. Therefore $\mathcal{M}$ is $i$-open star refinement of $\mathcal{U}$.\\

\begin{corollary} Let $X$ be pairwise almost regular, $B_1$-nearly paralindelof, pairwise $P_r$ and $P$-space. Then, every pairwise regular open covering of $X$ has open star refinement.\end{corollary}

\begin{remark} We  can notice that every RR-pairwise paralindelof is $B_1$-nearly paralindelof. But the converse is not true in general because not each $i$-open cover is $(i,j)$-regular open cover. \end{remark}

\begin{definition} \cite{Zabidin1}\cite{Zabidin6}. A bitopological space $X$ is said to be $(i,j)$-almost Lindelof (resp. $(i,j)$-weakly Lindelof) if for every $i$-open cover $\{U_\alpha :\alpha \in \Delta\}$ of $X$, there exists a countable subset $\{\alpha_n:n \in \mathbb{N}\}$ such that $X=\bigcup_{n \in \mathbb{N}} j-cl(U_{\alpha_n})$ (resp. $X=j-cl(\bigcup_{n \in \mathbb{N}} U_{\alpha_n})$). $X$ is called pairwise almost Lindelof (resp. pairwise weakly Lindelof) if it is both $(1,2)$-almost Lindelof and $(2,1)$-almost Lindelof (resp. $(1,2)$-weakly Lindelof and $(2,1)$-weakly Lindelof).\end{definition}

\begin{proposition} \label{pro1}. In $(i,j)_1$-nearly paralindelof,  $(i,j)$-semiregular and $j$-$P$-space, every $(i,j)$-weakly Lindelof is $(i,j)$-almost Lindelof. \end{proposition}

\textbf{Proof.} Let $\mathcal{U}=\{U_{\alpha}:\alpha \in \Delta\}$ be an $(i,j)$-regular open cover of $X$. Since $X$ is $(i,j)_1$-nearly paralindelof, $\mathcal{U}$ admits an $i$-open refinement which is $j$-locally countable $\mathcal{V}=\{V_\lambda:\lambda \in \Lambda\}$. Since $X$ is $(i,j)$-weakly Lindelof, there exists a countable subfamily $\{V_{\lambda_n}:n \in \mathbb{N}\}$ such that $X=j-cl(\cup_{n \in \mathbb{N}} V_{\lambda_n})$. Since $X$ is $j$-$P$-space and the family $\mathcal{V}*=\{V_{\lambda_n}:n \in \mathcal{N}\}$ is $j$-locally countable (since $\mathcal{V}* \subset \mathcal{V}$), then $j-cl(\bigcup_{n \in \mathbb{N}} V_{\lambda_n}) = \bigcup_{n \in \mathbb{N}} j-cl(V_{\lambda_n})$. Picking for $n \in \mathbb{N}, \alpha_n \in \Delta$ such that $V_{\lambda_n} \subseteq U_{\alpha_n}$, we get $X=\bigcup_{n \in \mathbb{N}} j-cl(V_{\lambda_n}) = \bigcup_{n \in \mathbb{N}} j-cl(U_{\alpha_n})$. Therefore, $X$ is $(i,j)$-almost Lindelof.\\

\begin{corollary} In $B_1$-nearly paralindelof, pairwise semiregular and $P$-space, every pairwise weakly Lindelof is pairwise almost Lindelof. \end{corollary}

\begin{proposition} (see \cite{Zabidin6}).\label{pro2} An $(i,j)$-regular space is $(i,j)$-almost Lindelof if and only if it is $i$-Lindelof.\end{proposition}

\begin{corollary} A pairwise regular space is pairwise almost Lindelof if and only if it is Lindelof. \end{corollary}

\begin{proposition} In $(i,j)_1$-nearly paralindelof, $(i,j)$-regular and $j$-$P$-space, every $(i,j)$-weakly Lindelof is $i$-Lindelof.\end{proposition}

\textbf{Proof.} Let $\{U_{\alpha}:\alpha \in \Delta\}$ be an $i$-open cover of $X$. Since $X$ is $(i,j)$-regular, then it is $(i,j)$-semiregular. So $X$ is $(i,j)$-almost Lindelof by Proposition \ref{pro1}. Also, since $X$ is $(i,j)$-regular, then $X$ is $i$-Lindelof by Proposition \ref{pro2}.\\

\begin{corollary} In $B_1$-nearly paralindelof, pairwise regular and $P$-space, every pairwise weakly Lindelof is Lindelof.\end{corollary}

\begin{definition} A bitopological space $(X,\tau_1,\tau_2)$ is said to be $(i,j)^*$-almost regular if for any $(i,j)$-regular closed  set $A$ and any point $x \in X$, there exist $i$-open set $U$ and $j$-open set such that $A \subseteq U, x \in V$ and $U \cap V =\emptyset$. $X$ is called $B^*$-almost regular if it is both $(1,2)^*$-almost regular and $(2,1)^*$-almost regular. \end{definition}

\begin{proposition} A bitopological space $X$ is $(i,j)^*$-almost regular space if and only if for every point $x \in X$ and for every $(i,j)$-regular open set $V$ of $X$ containing $x$, there exists an $(j,i)$-regular open set such that $x \in \subseteq U \subseteq i-cl(U) \subseteq V$.\end{proposition}

\textbf{Proof.} Let $V$ be any $(i,j)$-regular open set of $X$ containing $x$. Then, $X-V$ is $(i,j)$-regular closed set disjoint from $\{x\}$. Since $X$ is $(i,j)^*$-almost regular, there exist $i$-open set $A$ and $j$-open set $B$ in $X$ such that $x \in B, X- V \subset A$ and $A \cap B=\emptyset$. So $i-cl(B) \cap A=\emptyset \Longrightarrow i-cl(B) \subseteq X - A \subseteq V$. Then, $x \in B \subseteq i-cl(B) \subseteq V$. Again, $B \subseteq j-int(i-cl(B)) \subseteq i-cl(B) \subseteq V$. Set $U= j-int(i-cl(B))$, therefore, $x \in U \subseteq i-cl(U) \subseteq V$ as desire. Conversely, let $F$ be any $(i,j)$-regular closed set in $X$ and let any singleton $\{x\}$ disjoint from $F$. Hence, $X- F$ is $(i,j)$-regular open set containing $x$. By hypothesis, there exists $(j,i)$-regular open set $U$ in $X$ containing $x$ such that $x \in U \subseteq i-cl(U) \subseteq X\ F$. So $F \subseteq X- i-cl(U), x \in U$ and $U \cap X- i-cl(U)=\emptyset$. Therefore, $X$ is $(i,j)^{*}$-almost regular.\\

\begin{corollary} A bitopological space $X$ is $B^*$-almost regular space if and only if for every point $x \in X$ and for every pairwise regular open set $V$ of $X$ containing $x$, there exists an pairwise regular open set such that $x \in \subseteq U \subseteq i-cl(U) \subseteq V$ for each $i=1,2$.\end{corollary}

\begin{theorem} \label{almost regular} Let a bitopological space $X$ be $j$-Hausdorff and $j$-P-space. If $X$ is $(i,j)_1$-nearly paralindelof, then it is $(i,j)^*$-almost regular.\end{theorem}

\textbf{Proof.} Let $x \in X$ and $F$ be any $(i,j)$-regular closed subset of $X$ such that $ x \notin F$. Now, for every  $y \in F$, since $X$ is $j$-Hausdorff space, there exist two disjoint $j$-open sets $ U_{xy}$ and $V_y$ such that  $x \in U_{xy} , y \in V_y$. Since $ U_{xy} \cap V_y = \emptyset \Longrightarrow j-cl(V_y) \cap U_{xy} = \emptyset$, $y \notin U_{xy}$. Furthermore, $\{x\} \subseteq U_{xy}$ so that $\{x\} \cap V_y= \emptyset \Longrightarrow j-cl(V_y) \subseteq X \ \{x\}$.\\

\noindent Now $\mathcal{U}=\{ i-int(j-cl(V_y)):y \in F \} \bigcup \{X -F\}$ is $(i,j)$-regular open cover of $X$. Then $\mathcal{U}$ has $j$-locally countable family $\mathcal{W}= \{W_\alpha : \alpha \in \Delta \} \bigcup \{X - F\}$ of $i$-open sets which refines $ \mathcal{U}$. Set $V=\bigcup\{W_\alpha : \alpha \Delta\}$ so that $V$ is $i$-open set containing $F$. Further, let $H=X -\bigcup \{j-cl(W_\alpha : \alpha \in \Delta\}$. Since $\{W_\alpha : \alpha \in \Delta\}$ is $j$-locally countable and $X$ is $j$-$P$-space, then $\bigcup \{j-cl(W_\alpha): \alpha \in \Delta\}$ is $j$-locally countable by Theorem \ref{locally countable}. Thus, $H$ is $j$-open set. \\

Since $\mathcal{W}$ refines $\mathcal{U}$ and every $W_\alpha$ meets $F$, for each $\alpha \in \Delta$, there is $y \in F$ such that $W_\alpha \subseteq i-int(j-cl(V_y)) \Longrightarrow j-cl(W) \subseteq j-cl(i-int(j-cl(V_y))) \subseteq j-cl(V_y) \subseteq X -\{x\}$. Hence, $x \notin j-cl(W)$ for each $\alpha \in \Delta$. Thus, $x \notin \bigcup \{j-cl(W_\alpha : \alpha \in \Delta\}$. So $x \in X \ \bigcup \{j-cl(W_\alpha : \alpha \in \Delta\}=H$. Then, $F \subseteq V, x \in H$ $V \cap H=\emptyset$. Therefore, $X$ is $(i,j)^*$-almost regular.\\

\begin{corollary} Let a bitopological space $X$ be Hausdorff and P-space. If $X$ is $B_1$-nearly paralindelof, then it is $B^*$-almost regular.\end{corollary}

\begin{definition} A bitopological space $X$ is said to be $(i,j)$-weakly regular if for each $(i,j)$-regular closed set $F$ and for each  $x \notin F$ such that $j-cl\{x\} \cap F=\emptyset$, there exist $i$-open set $U$ and $j$-open set such that $\{x\} \subseteq U, F \subseteq V$ and $U \cap V=\emptyset$. $X$ is called pairwise weakly regular if it is both $(1,2)$-weakly regular and $(2,1)$-weakly regular. \end{definition}

\begin{theorem} A bitopological space $X$ is $(i,j)$-weakly regular if and only if for each point $x \in X$ and for each $(i,j)$-regular open set $U$ containing $j-cl\{x\}$, there exists $i$-open set $V$ such that $x \in V \subseteq j-cl(V) \subseteq U$. \end{theorem}

\textbf{Proof.} Suppose that $X$ is $(i,j)$-weakly regular. Let $x \in X$ and $V$ be any $(i,j)$-regular open set containing $j-cl\{x\}$. So, $G=X -V$ is $(i,j)$-regular closed set and $x \notin G$. Since $X$ is $(i,j)$-weakly regular, there exist there exist $i$-open set $A$ and $j$-open set $B$ such that $\{x\} \subseteq A, G \subseteq B$ and $A \cap B=\emptyset$. Now, $A \subseteq X-B \Longrightarrow j-cl(A) \subseteq j-cl(X -B)=X -B \subset X -G=V$. Then, we have, $x \in A \subseteq j-cl(A) \subseteq V$. Hence the result. Conversely, let the condition holds. Let $x \in X$ and $F$ be any $(i,j)$-regular closed such that $x \notin F$. Thus, $x \in X \ F$ and by hypothesis, there exists $i$-open set $U$ such that $x \in U \subseteq j-cl(U) \subseteq X -F$. Then, $x \in U, F \subseteq X-j$-$cl(U)$ and $U \cap X -j$-$cl(U)=\emptyset$. Therefore, $X$ is $(i,j)$-weakly regular.\\

\begin{corollary} A bitopological space $X$ is pairwise weakly regular if and only if for each point $x \in X$ and for each pairwise regular open set $U$ containing $j-cl\{x\}$, there exists $i$-open set $V$ such that $x \in V \subseteq j-cl(V) \subseteq U$ for each $i,j=1,2$ and $i \neq j$. \end{corollary}

\begin{remark} \label{weakly+almost} Clearly, every $(i,j)$-almost regular space is $(i,j)$-weakly regular because each $(i,j)$-regular open set is $i$-open. Furthermore, if the bitopological space $X$ is $T_1$, then, the concept of pairwise weakly regularity and pairwise almost regularity are equivalent. \end{remark}

\begin{theorem} In an $(i,j)$-weakly regular, $(i,j)_1$-nearly paralindelof and $j$-P-space, every pair disjoint $(i,j)$-regular closed sets can be strongly separated.\end{theorem}

\textbf{Proof.} Let $F$ and $G$ be any $(i,j)$-regular closed subsets of $X$. Then, for every point $x \in F$, such that $j-cl(\{x\})$ contained in $F$. Since $X$ is $(i,j)$-weakly regular, there exists $i$-open set $U_x$ containing $x$ and $j-cl(U_x) \cap G=\emptyset$. Now, $\mathcal{U}=\{U_x: x \in F\} \bigcup \{X-F\}$ is a family of $(i,j)$-regular open sets which covers $X$. Hence, it has $j$-locally countable family $\mathcal{W}=\{W_\alpha : \alpha \in \Delta\} \bigcup \{X-F\}$ of $i$-open sets refining $\mathcal{U}$ and covering $X$. Let $W=\bigcup \{W_\alpha : \alpha \in \Delta\}$ so that $W$ is $i$-open set containing $F$. Set $H=X-\bigcup \{j-cl(W_\alpha) : \alpha \in \Delta\}$. As the proof of Theorem \ref{almost regular} and by Theorem \ref{locally countable}, $H$ is $j$-open set. For each $\alpha \in \Delta$, there exists $x \in F$ such that $W_\alpha \subseteq i-int(j-cl(U_x)) \Longrightarrow l-cl(W_\alpha) \subseteq j-cl(i-int(j-cl(U_x))) \subseteq j-cl(U_x) \subseteq X-G$. Then, $G \cap j-cl(W_\alpha)=\emptyset$ for all $\alpha \in \Delta$. Hence, $G \subseteq X-\bigcup \{j-cl(W_\alpha): \alpha \in \Delta\}=H$. Therefore, $W$ and $H$ are $i$-open and $j$-open set respectively such that $F \subseteq W, G \subseteq H$ and $W \cap H=\emptyset$ and hence the result.

\begin{corollary} In pairwise weakly regular, $B_1$-nearly paralindelof and P-space, every pair disjoint pairwise regular closed sets can be strongly separated.\end{corollary}

\begin{corollary} In an $(i,j)$-almost regular, $(i,j)_1$-nearly paralindelof and $j$-P-space, every pair disjoint $(i,j)$-regular closed sets can be strongly separated.\end{corollary}

\textbf{Proof.} From Remark \ref{weakly+almost}, every $(i,j)$-almost regular space is $(i,j)$-weakly regular.

\begin{corollary} In an pairwise almost regular, $B_1$-nearly paralindelof and P-space, every pair disjoint pairwise regular closed sets can be strongly separated.\end{corollary}

\begin{definition} \label{Br} A bitopological space $(X,\tau_1,\tau_2)$ is said to be $(i,j)_r$-nearly paralindelof if every $(i,j)$-regular open cover $\{U_\alpha:\alpha \in \Delta\}$ of $X$ admits $(i,j)$-regular open $j$-locally countable refinement $\{V_\beta:\beta \in B\}$ such that $j-cl(V_\beta) \subseteq U_{\alpha(\beta)}$ for some $\beta \in B$. $X$ is called $B_r$-nearly paralindelof is it is both $(1,2)_r$-nearly paralindelof and $(2,1)_r$-nearly paralindelof.\end{definition}

\begin{remark} It clear that every $(i,j)_r$-nearly paralindelof space is $(i,j)_1$-nearly paralindelof since every $(i,j)$-regular open set is $i$-open. But the converse is not true in general.\end{remark}

\begin{theorem} Let $X$ be $(i,j)$-almost regular $(i,j)$-semiregular and $j$-$P$-space. Then, $X$ is $\tau_i$ paralindelof respect to $\tau_j$ if and only if $X$ is $(i,j)_r$-nearly paralindelof.\end{theorem}

\textbf{Proof.} Let $X$ be $\tau_i$ paralindelof respect to $\tau_2$. Let $\mathcal{U}=\{U_\alpha :\alpha \in \Delta\}$ be $(i,j)$-regular open cover of $X$. For $x \in X$, there is a set $U_\alpha \in \mathcal{U}$ containing $x$. Since $X$ is $(i,j)$-almost regular, there exists $(i,j)$-regular open set neighbourhood $V_x$ of $x$ such that $x \in V_x \subseteq j-cl(V_x) \subseteq U_\alpha$. The collection $\mathcal{V}=\{V_x: x \in X\}$ is $(i,j)$-regular open cover of $X$ so that it is also $i$-open cover of $X$. By hypothesis, $\mathcal{V}$ has $j$-locally countable family $\mathcal{W}$ of $i$-open sets refining $\mathcal{V}$ and covering $X$, i.e., $X = \bigcup \{W : W \in \mathcal{W}\}  \subseteq \bigcup \{i-int(j-cl(W))  W \in \mathcal{W}\}$ so that $\mathcal{H}=\{i-int(j-cl(W))=H: W \in \mathcal{W}\}$ is a family of $(i,j)$-regular open subsets which covers $X$. Now, we have to show that $\mathcal{H}$ is $j$-locally countable, refines $\mathcal{U}$ and $j-cl(H) \subseteq U_\alpha$ for some $\alpha \in \Delta$. Since $ \{i-int(j-cl(W)):W \in \mathcal{W}\} \subseteq \{j-cl(W): W \in \mathcal{W}\}$ and $\{j-cl(W): W \in \mathcal{W}\}$ is $j$-locally countable, $\mathcal{H}$ is also $j$-locally countable by Theorem \ref{locally countable}. For each $W \in \mathcal{W}$, pick some $\alpha_W \in \Delta$ such that $H \subseteq j-cl(H) \subseteq j-cl(W) \subseteq U_{\alpha_W}$. Therefore, $X$ is $(i,j)_r$-nearly paralindelof.\\
Conversely, in $(i,j)$-semiregular space, $(i,j)$-regular open set and $i$-open are equivalent. Thus, the proof completes.

\begin{corollary} Let $X$ be pairwise almost regular pairwise semiregular and $P$-space. Then, $X$ is $RR$-pairwise paralindelof if and only if $X$ is $B_r$-nearly paralindelof.\end{corollary}

\begin{lemma} \label{Br} Let $X$ be $(i,j)_r$-nearly paralindelof, $(i,j)$-almost regular and $j$-$P$-space. Then, every $(i,j)$-regular open cover $\mathcal{U}=\{U_\alpha: \alpha \in \Delta\}$ of $X$ has $(i,j)$-regular open refinement $j$-locally countable family $\mathcal{V}=\{V_\alpha: \alpha \in \Delta\}$ such that $j-cl(V_\alpha) \subseteq U_{\alpha}$ for all $\alpha \in \Delta$.\end{lemma}

\textbf{Proof.} Since $X$ is $(i,j)$-almost regular space, there is $(i,j)$-regular open cover $\mathcal{W}=\{W: W
\in \mathcal{W}\}$ refining $\mathcal{U}$ such that, for all $W \in \mathcal{W}$, $j-cl(W) \subseteq U_{\alpha}$ for some $\alpha \in \Delta$. Then, $\mathcal{W}$ has $j$-locally countable family $\mathcal{G}=\{G_\lambda:\lambda \in \Lambda\}$ of $(i,j)$-regular open sets which refines $\mathcal{W}$ and covers $X$ such that for all $\lambda \in \Lambda$ we have $j-cl(G_\lambda) \subseteq U_{\alpha}$ for some $\alpha \in \Delta$. Take $f:\Lambda \rightarrow \Delta$. Now, for each $\lambda \in \Lambda$, choose $f(\lambda) \in \Delta$ such that $j-cl(G_\lambda) \subseteq U_{f(\lambda)}$. Set $\Lambda_{\alpha}=f^{-1}(\alpha)$ for all $\alpha \in \Delta$.\\

Notice that some members of $\Lambda_{\alpha}$ might be empty and $\Lambda_{\alpha} \cap \Lambda_{\alpha^{'}}=\emptyset$ for $\alpha \neq \alpha^{'}$. Let $H_\alpha=\bigcup_{\lambda \in \Lambda_{\alpha}} G_\lambda$. So, the family $\mathcal{H}=\{ H_\alpha:\alpha \in \Delta \}$ covers $X$. Now, by using Theorem \ref{closure preserving}, we have
$j-cl(H_\alpha)= j-cl(\bigcup_{\lambda \in \Lambda_{\alpha}} G_\lambda)=\bigcup_{\lambda \in \Lambda_{\alpha}} j-cl(G_\lambda) \subseteq U_\alpha$

hence $j-cl(H_\alpha) \subseteq U_\alpha$ for all $\alpha \in \Delta$. Then, $\mathcal{V}=\{i-int(j-cl(H_\alpha))=V_\alpha: \alpha \in \Delta\}$ is $(i,j)$-regular open cover of $X$ since $H_\alpha \subseteq i-int(j-cl(H_\alpha))$. Since $\mathcal{V}$ is a subfamily of $j$-locally countable family, then $\mathcal{V}$ is $j$-locally countable such that $j-cl(V_\alpha)=j-cl(i-int(j-cl(H_\alpha))) \subseteq j-cl(j-cl(H_\alpha))=j-cl(H_\alpha) \subseteq U_\alpha$ for each $\alpha \in \Delta$.

\begin{corollary} Let $X$ be $B_r$-nearly paralindelof, pairwise almost regular and $j$-$P$-space. Then, every pairwise regular open cover $\mathcal{U}=\{U_\alpha: \alpha \in \Delta\}$ of $X$ has pairwise regular open refinement $j$-locally countable family $\mathcal{V}=\{V_\alpha: \alpha \in \Delta\}$ such that $j-cl(V_\alpha) \subseteq U_{\alpha}$ for all $\alpha \in \Delta$ for all $i,j=1.2$ and $i \neq j$.\end{corollary}

\begin{theorem} Let a bitopological space $X$ be $(i,j)$-almost regular, $(j,i)$-weakly $P$ and $j$-$P$-space. If $X$ is $(i,j)_r$-nearly paralindelof such that there exists $j$-dense set $A$ of $X$ which is $(i,j)$-nearly Lindelof subspace, then $X$ is $(i,j)$-nearly Lindelof space.\end{theorem}

\textbf{Proof.} Let $\mathcal{U}=\{U_\alpha: \alpha \in \Delta\}$ be any $(i,j)$-regular open cover of $X$. By Lemma \ref{Br}, there exists $(i,j)$-regular open $j$-locally countable refinement $\mathcal{V}=\{V_\alpha:\alpha \in \Delta\}$ of $\mathcal{U}$ such that $j-cl(V_\alpha) \subseteq U_\alpha$ for each $\alpha \in \Delta$. Then, $\{A \cap V_\alpha: \alpha \in \Delta\}$ forms $(i,j)$-regular open cover of subspace $A$. Therefore, there exists a countable set $\Delta_0 \subset \Delta$ such that $A \subseteq \bigcup\{A \cap V_\alpha: \alpha \in \Delta_0\}$ Then, since $X$ is $ij$-weakly $P$-space, we have
\begin{center}
$X=j-cl(A)=j-cl(\bigcup_{\alpha \in \Delta_0} (A \cap V_\alpha))= \bigcup_{\alpha \in \Delta_0} j-cl(A \cap V_\alpha) \subseteq \bigcup_{\alpha \in \Delta_0} j-cl(V_\alpha) \subseteq \bigcup_{\alpha \in \Delta_0} U_\alpha$,

\end{center}
Thus, $X$ is $(i,j)$-nearly Lindelof space.

\begin{corollary} Let a bitopological space $X$ be pairwise weakly $P$-space. If $X$ is $B_r$-nearly paralindelof such that there exists dense set $A$ of $X$ which is pairwise nearly Lindelof subspace, then $X$ is pairwise nearly Lindelof space.\end{corollary}

\section{pairwise almost paralindel$\ddot{o}$f subset}

In this section, we are going to study a new class of pairwise almost paralindelof namely $B_1$-almost paralindelof space.

\begin{definition} A bitopological space $X$ is $(i,j)_1$-almost paralindelof if every $i$-open cover $\mathcal{U}$ of $X$ admits $j$-locally countable collection $\mathcal{V}$ of $i$-open subsets of $X$ refining $\mathcal{U}$ such that  $X=\bigcup \{j-cl(V):V \in \mathcal{V}\}$. $X$ is called $B_1$-almost paralindelof if it is $(1,2)_1$-almost paralindelof and $(2,1)_1$-almost paralindelof.\end{definition}

\begin{proposition} Let $X$ be $(i,j)$-almost regular, $(i,j)_1$-almost paralindelof and $j$-$P$-space. Then, every $(i,j)$-regular open cover $\mathcal{U}$ of $X$, there exists $j$-locally countable family of $(j,i)$-regular closed sets which refines $\mathcal{U}$ and covers $X$.\end{proposition}

\textbf{Proof.} Let $\mathcal{U}$ be any  $(i,j)$-regular open cover of $X$. For each point $x \in X$, there exists $(i,j)$-regular open set $V_x$ such that
\begin{center}
$x \in V_x \subseteq j-cl(V_x) \subseteq U$
\end{center}
for some $U \in \mathcal{U}$. Now, $\mathcal{V}=\{V_x :x \in X\}$ is $(i,j)$-regular open cover of $X$. There exists $j$-locally countable family $\mathcal{W}$ of $i$-open sets which refines $\mathcal{V}$  such that $X = \bigcup \{j-cl(W): W \in \mathcal{W}\}$. Now, for each $\alpha \in \Delta$, there exists $x(\alpha) \in X$ such that
\begin{center}
$W_\alpha \subseteq j-cl(V_{x(\alpha)}) \subseteq U$
\end{center}
for some $U \in \mathcal{U}$. Hence, we obtain

\begin{center}
$j-cl(W_\alpha) \subseteq j-cl(V_{x(\alpha)}) \subseteq U$.
\end{center}

Since $X$ is $j$-$P$-space, then the family $\{j-cl(W_\alpha): \alpha \in \Delta\}$ is $j$-locally countable collection, using Theorem \ref{locally countable}, of $(j,i)$-regular closed sets which refines $\mathcal{U}$ and covers $X$.\\

\begin{corollary} Let $X$ be pairwise almost regular, $B_1$-almost paralindelof and $P$-space. Then, every pairwise regular open cover $\mathcal{U}$ of $X$, there exists locally countable family of pairwise regular closed sets which refines $\mathcal{U}$ and covers $X$.\end{corollary}

\begin{remark} It is clear that every $RR$-pairwise paralindelof spaces are $B_1$-almost paralindelof.\end{remark}

\begin{theorem} Let $X$ be $(i,j)$-almost regular, $(j,i)$-extremely disconnected and $j$-$P$-space. Then, $X$ is $(i,j)_1$-nearly paralindelof if it is $(i,j)_1$-almost paralindelof. \end{theorem}

\textbf{Proof.} If $X$ is $(i,j)_1$-nearly paralindelof, then clearly that it is $(i,j)_1$-almost paralindelof. Conversely, let $\mathcal{U}=\{U_\alpha : \alpha \in \Delta\}$ be any  $(i,j)$-regular open cover of $X$. For each $x \in X$, there is $(i,j)$-regular open set $V_x$ such that
\begin{center}
$x \in V_x \subseteq j-cl(V_x) \subseteq U_\alpha$
\end{center}
for some $\alpha \in \Delta$. Now, $\mathcal{V}=\{V_x: x \in X\}$ forms $(i,j)$-regular open cover of $X$. Since $X$ is $(i,j)_1$-almost paralindelof, there exists $j$-locally countable family $\mathcal{W}=\{W_\beta: \beta \in B\}$ of $i$-open sets that refines $\mathcal{V}$ such that $X=\bigcup \{j-cl(W_\beta): \beta \in B\}$. Set $\mathcal{H}=\{j-cl(W_\beta): \beta \in B\}$. For every $\beta \in B$, there exists $x(\beta) \in X$ such that $W_\beta \subseteq j-cl(W_\beta) \subseteq j-cl(V_x) \subseteq U_\alpha(x)$ for some $\alpha(x) \in \Delta$. Hence $\{j-cl(W_\beta): \beta \in B\}$ refines $\mathcal{U}$. Since $X$ is $(j,i)$-extremely disconnected, then $j-cl(W_\beta)=i-int(j-cl(W_\beta))$, i.e., $j-cl(W_\beta)$ is $i$-open set, say $H$. Moreover, $\{j-cl(W_\beta): \beta \in B\}$ is $j$-locally countable because $X$ is $j$-$P$-space. Therefore, $\mathcal{H}$ is $j$-locally countable family of $i$-open sets refining $\mathcal{U}$. Thus, $X$ is $(i,j)_1$-nearly paralindelof.\\

\begin{corollary} Let $X$ be pairwise almost regular, pairwise extremely disconnected and $P$-space. Then, $X$ is $B_1$-nearly paralindelof if it is $B_1$-almost paralindelof.\end{corollary}

\end{document}